\documentclass[12pt]{article}

\usepackage{amssymb}
\usepackage{amsmath}

\begin{document}%



\title{\bf Ramanujan's Master Theorem and two formulas for zero-order Hankel transform}

\author{
A.V. Kisselev\thanks{Electronic address:
alexandre.kisselev@ihep.ru} \\
{\small A.A.~Logunov Institute for High Energy Physics, NRC
``Kurchatov
Institute''}, \\
{\small 142281 Protvino, Russia}
}

\date{}

\maketitle

\begin{abstract}
Using Ramanujan's Master Theorem, two formulas are derived which
define the Hankel transforms of order zero with even functions by
inverse Mellin transforms, provided these functions and their
derivatives obey special conditions. Their validity is illustrated
by a number of examples. With a help of one of these formulas, one
as yet unknown parametric improper integral of the Bessel function
$J_0(x)$ is calculated.
\end{abstract}



\section{Introduction} %
\label{sec:int}

The Hankel transform,
\begin{equation}\label{Hankel_tramsform}
\mathcal{H}_\nu = \int\limits_0^\infty \!\!J_\nu (q x) f(x) \,x \,dx
\;,
\end{equation}
where $J_\nu (x)$ is the Bessel function of order $\nu$, solves a
number of problems in mathematical physics \cite{Tranter},
\cite{Sneddon}, \cite{Debnath} and high energy nuclear and particle
physics \cite{Moliere}, \cite{Glauber}, \cite{Collins}.

Often it is necessary to know an asymptotics of this integral as $q
\rightarrow \infty$. A method of obtaining such an asymptotics were
done in \cite{Willis} (see also \cite{Tranter}), where the following
results were obtained:
\begin{align}
\int\limits_0^\infty \!\!J_0(q x) f(x) \,dx &= \frac{f(0)}{q} -
\frac{1}{2} \frac{f^{(2)}(0)}{q^3} + \frac{1 \cdot 3}{2^2 \cdot 2!}
\frac{f^{(4)}(0)}{q^5}  \nonumber \\
&- \frac{1 \cdot 3 \cdot 5}{2^3 \cdot 3!}
\frac{f^{(6)}(0)}{q^7} + \ldots \;, \label{Willis_J0} \\
\int\limits_0^\infty \!\!J_1(q x) f(x) \,dx &= \frac{f(0)}{q} +
\frac{f^{(1)}(0)}{q^2} - \frac{1}{2} \frac{f^{(3)}(0)}{q^4} +
\frac{1 \cdot 3}{2^2 \cdot 2!} \frac{f^{(5)}(0)}{q^6} + \ldots
\label{Willis_J1}\;.
\end{align}
Note that integrands in integrals \eqref{Willis_J0},
\eqref{Willis_J1}, in contrast to \eqref{Hankel_tramsform}, do not
contain a factor $x$. Later on, asymptotic expansions for
$\mathcal{H}_\nu$ as $q \rightarrow \infty$ have been obtained by
several authors \cite{Hsu}-\cite{Soni}.

It follows from \eqref{Willis_J0}  that for $\nu=0$ Hankel transform
\eqref{Hankel_tramsform} is presented by an asymptotic series which
contains derivatives of \emph{odd orders} only:
\begin{equation}\label{asym_J0}
\mathcal{H}_0(q) = \frac{1}{q^3} \sum_{m=0}^\infty (-1)^{m+1}
\,\frac{\Gamma (2m+2)}{\Gamma^2 (m+1)} f^{(2m+1)}(0) \, (2q)^{-2m}
\;.
\end{equation}
It is clear that a case when $f^{(2m+1)}(0) = 0$ for \emph{all}
$m=0,1,2, \ldots$ should be considered separately, as was mentioned
in \cite{Frenzen}. Then one expects that the Hankel transform
$\mathcal{H}_0(q)$ decreases as $q \rightarrow \infty$ more quickly
than any inverse power of $q$. For instance, it could take place if
$f(x)$ is an \emph{even} function, $f(x) = g(x^2)$.

The main goal of this paper is to derive two representations for the
Hankel transform of order zero with even functions which are valid
not only for $q \rightarrow \infty$ but for \emph{all} $q>q_0>0$. As
a byproduct, one as yet unknown parametric improper integral of
$J_0(x)$ will be calculated.

\section{Hankel transform of order zero} %
\label{sec:J0}

Let us study the Hankel transform of order zero with an \emph{even}
function:
\begin{equation}\label{A}
A(q) = \mathcal{H}_0(q) = \int\limits_0^\infty \!\!J_0(q x) f(x) \,x
dx \;,
\end{equation}
where $f(x) = g(x^2)$.%
\footnote{A special case $f(x)= e^{-x^2}\varphi(x^2)$ was considered
in \cite{Frenzen}.}

\textsc{Theorem 1.} The Hankel transform of order zero \eqref{A}
with the function $f(x) = g(x^2)$ and $q>0$ may be expressed in the
form
\begin{equation}\label{A_cont_final}
A(q) = \frac{1}{\pi i q^2}\int\limits_{\alpha - i\infty}^{\alpha +
i\infty} \!\!ds \,\bar{g}^{(s)}(0) \,\Gamma (s+1) \!\left(
\frac{q^2}{4} \right)^{\!\!-s} ,
\end{equation}
for $-1 < \alpha < 0$, provided that
\begin{description}
  \item[(1)]
$g(z)$ is a regular function and its Taylor series at $z=0$  has the
form
\begin{equation}\label{g_Taylor_series}
g(z) = \sum_{m=0}^\infty \frac{\bar{g}^{(m)}(0)}{m!} \,(-z)^m \;;
\end{equation}
  \item[(2)]
$g(z)= \mathrm{O}(z^{-d})$, as $z\rightarrow \infty$, for $d>1/4$;
  \item[(3)]
$\bar{g}^{(s)}(0)$ is a regular (single-valued) function defined on
a half-plane
\begin{equation}\label{g_half_plane}
H(\delta) = \{s\in \mathrm{\mathbb{C}}\!: \mathrm{Re} \,s \geqslant
-\delta \}
\end{equation}
for some $1/4 < \delta < 1$ and satisfies the growth condition
\begin{equation}\label{g_growth}
|\bar{g}^{(s)}(0)| < C e^{P v + A|w|}
\end{equation}
for some $A < \pi/2$ and all $s=v+iw \in H(\delta)$.
\end{description}

\textsc{Proof}. The main tool of the proof will be the power
\emph{Ramanujan Master Theorem} (see, for instance, its proof
provided by Hardy in \cite{Hardy}). It provides an analytic
expression for the Mellin transform of an analytic function. We will
use the following version of Ramanujan's master theorem
\cite{Amdeberham}:
\begin{equation}\label{master_theorem}
\int\limits_0^\infty x^{s-1} \left[ \varphi(0) - x \varphi(1) + x^2
\varphi(2) - \ldots \right] dx = \frac{\pi}{\sin (s\pi)} \, \varphi
(-s) \;,
\end{equation}
where $0 < \mathrm{Re} \,s < \delta$. Formula \eqref{master_theorem}
is valid if the function $\varphi (s)$ in \eqref{master_theorem}
satisfies condition~3 of the Theorem~1 with $A < \pi$.

Let us put
\begin{equation}\label{phi}
\varphi (s) = \frac{\bar{g}^{(s)}(0)}{\Gamma(s+1)}
\end{equation}
in \eqref{master_theorem}. One can see that this function obeys all
conditions of Ramanujan's Master Theorem, and we have
\begin{equation}\label{master_theorem_appl}
\int\limits_0^\infty x^{s-1} \left[ \sum_{m=0}^\infty
\frac{\bar{g}^{(m)}(0)}{\Gamma(m+1)} \,(-x)^m  \right] =
\frac{\pi}{\sin (s\pi)} \, \frac{\bar{g}^{(-s)}(0)}{\Gamma(-s+1)}
\;.
\end{equation}
Then the inverse Mellin transform gives immediately ($1/4 < c <
\delta$)
\begin{align}\label{Mellin_transform}
g(z) &= \frac{1}{2i} \int\limits_{c-i\infty}^{c+i\infty}
\frac{ds}{\sin(s\pi)} \frac{\bar{g}^{(-s)}(0)}{\Gamma(-s+1)}
\,z^{-s} \nonumber \\
&= -\frac{1}{2i} \int\limits_{-c-i\infty}^{-c+i\infty}
\frac{ds}{\sin(s\pi)} \frac{\bar{g}^{(s)}(0)}{\Gamma(s+1)} \,z^s \;.
\end{align}

Now we put $z=x^2$ and replace $g(x^2)$ in \eqref{A} by its integral
representation \eqref{Mellin_transform}. As a result, we obtain the
formula
\begin{equation}\label{A_Mellin_Bessel}
A(q) = -\frac{1}{2i} \int\limits_{-c-i\infty}^{-c+i\infty}
\frac{ds}{\sin(s\pi)} \frac{\bar{g}^{(s)}(0)}{\Gamma(s+1)}
\int\limits_0^\infty dx \,x^{2s+1} \, J_0(q x) \;.
\end{equation}
Using equation \cite{Bateman_vol_2}
\begin{equation}\label{int_Bessel}
\int\limits_0^\infty  dx \,x^{2s+1} J_0(q x) = 2^{\,2s+1} \,
q^{-2s-2} \, \frac{\Gamma(1+s)}{\Gamma(-s)} \;,
\end{equation}
valid for $-1< - \delta < \mathfrak{Re} s < -1/4$, as well as
equation \cite{Bateman_vol_1}
\begin{equation}\label{Gamma_relation}
\frac{1}{\sin(s \pi) \Gamma(-s)} = - \frac{1}{\pi} \,\Gamma(s+1) \;,
\end{equation}
we come to formula \eqref{A_cont_final}. Q.E.D.

Let us see how does our formula \eqref{A_cont_final} work? It is
worth it to consider a number of examples (everywhere below
it is assumed that $q>0$). \\
\textbf{A.1.} $f(z) = e^{-a^2 x^2}$, $a>0$, then $\bar{g}^{(m)}(0) =
a^{2m}$, and we get from \eqref{A_cont_final} ($\alpha
> -1$)
\begin{align}\label{example_cont_int_A1}
A_1(q) &= \int\limits_0^{\infty} \!\!x e^{-a^2 x^2} \! J_0(q x) \,dx
\nonumber \\
&= \frac{1}{\pi i q^2}\int\limits_{\alpha - i\infty}^{\alpha +
i\infty} \!\!ds \,\Gamma(1+s) \,\left( \frac{q^2}{4a^2}
\right)^{\!\!-s} = \frac{1}{2a^2} \,e^{-q^2/4a^2} \;,
\end{align}
in accordance with eq.~2.12.9.3. in \cite{Prudnikov_vol_2}. To get
this result, we used formula 7.3(1) from
ref.~\cite{Bateman_tab_vol_1} ($\mathrm{Re} \,c > 0$):
\begin{equation}\label{contour_int_A1}
\frac{1}{2\pi i}\int\limits_{c - i\infty}^{c + i\infty} \!\!ds \,
(ax)^{-s} \,\Gamma(s) = e^{-ax} \;.
\end{equation}
\textbf{A.2.} $f(x) = e^{-a^2 x^2} J_0(c x)$, $a, c>0$. The
derivatives $\bar{g}^{(m)}(0)$ are calculated in Appendix~A to be
\begin{equation}\label{deriv_A2}
\bar{g}^{(m)}(0) =  a^{2m} L_m \!\left( - \frac{c^2}{4a^2}\right)
\;,
\end{equation}
where $L_n(z)$ is the Laguerre polynomial \cite{Bateman_vol_2}. Note
that $L_{n}(-z) > 0$ for $n \geqslant 0$ if $z>0$. It means that
$\bar{g}^{(m)}(0) > 0$ for all $m\geqslant 0$.

We have the relation \cite{Bateman_vol_1}
\begin{equation}\label{laguerre_whitakker_0}
L_m(-z)= e^{z/2} (-z)^{-1/2} M_{m+1/2,\,0} (-z) \;,
\end{equation}
where $M_{\lambda,\chi}(z)$ is the Whittaker function of the first
kind \cite{Bateman_vol_1}. Thus, we obtain ($-1 < \alpha <0$)
\begin{align}\label{int_A2}
A_2(q) &= \int\limits_0^{\infty} \!\!x e^{-a^2 x^2} \!J_0(c x) J_0(q
x) \,dx = \left( \!- \frac{c^2}{4a^2}\right)^{-1/2}
\!\!e^{-c^2/8a^2}
\nonumber \\
&\times \frac{1}{\pi i q^2}\int\limits_{\alpha - i\infty}^{\alpha +
i\infty} \!\!ds \,\Gamma (s+1) M_{s+1/2,\,0} \!\left( \!-
\frac{c^2}{4a^2}\right) \!\left( \frac{q^2}{4a^2} \right)^{\!\!-s} .
\end{align}

Now we can apply formula 7.5(19) from \cite{Bateman_tab_vol_1}%
\footnote{The tabulated expression for this contour integral
presented in \cite{Bateman_tab_vol_1} has a misprint. Namely, the
factor $x^{1/2}$ is missing in it.}
\begin{align}\label{whittaker_con_int}
&\frac{1}{2\pi i}\!\int\limits_{c - i\infty}^{c + i\infty} \!\!ds \,
\Gamma (s+\nu+1/2) M_{s,\,\nu}(y) \,x^{-s}
\nonumber\\
&= \Gamma(2\nu+1) (xy)^{1/2} \, e^{-x + y/2} J_{2\nu}\!\left(
2(xy)^{1/2} \right) ,
\end{align}
valid for $\mathrm{Re}(c + \nu)> -1/2$. After replacement $s+1/2 =
s'$ in \eqref{int_A2} and taking into account that $J_0(iz) =
I_0(z)$, we find (see eq.~2.12.39.3. in \cite{Prudnikov_vol_2})
\begin{equation}\label{example_int_A2}
A_2(q)  =  \frac{1}{2a^2} \exp \!\left( - \frac{q^2+c^2}{4a^2}
\right) \!I_0 \!\left( \frac{q c}{4a^2} \right),
\end{equation}
where $I_\nu(z)$ is the modified Bessel function of the first kind.

As it turns out, our formula \eqref{A_cont_final} gives correct
results even if the condition~3 of the Theorem~1 is violated
(namely, if $A$ \emph{is equal to} $\pi/2$ in \eqref{g_growth}).

Let us illustrate this statement by several examples.\\
\textbf{A.3.} $f(x) = 1/(x^2+a^2)^{n+1}$, $a>0$, integer
$n\geqslant0$. Then
\begin{equation}\label{example_int_A3}
A_3(q) = \int\limits_0^{\infty} \!\!\frac{x}{(x^2 + a^2)^{n+1}}
\,J_0(q x) \,dx \;.
\end{equation}
We have $\bar{g}^{(s)}(0) = a^{-2m-2n-2} \,\Gamma(m+n+1)/\Gamma(n+1)$.%
\footnote{Note that $\Gamma(s+n+1)$ grows as $\exp(A \mathfrak{Im}
s)$ with $A=\pi/2$ when $\mathfrak{Im}s \rightarrow -\infty$ .}
Then we obtain from \eqref{A_cont_final} (for $-1 < \alpha <0$)
\begin{align}\label{A3}
A_3(q) &= \frac{1}{a^{2n+2} q^2} \frac{1}{\Gamma(n+1)} \frac{1}{\pi
i}\int\limits_{\alpha - i\infty}^{\alpha + i\infty} ds \,\Gamma(s+1)
\,\Gamma(s+n+1) \left( \frac{q a}{2}\right)^{\!-2s} \nonumber \\
&= \frac{1}{\Gamma(n+1)}  \left( \frac{q}{2a}\right)^{\!n} \!K_{n}(q
a) \;,
\end{align}
where $K_\nu(z)$ is the modified Bessel function of the second kind
(Macdonald function), in accordance with eq.~2.12.4.28. in
\cite{Prudnikov_vol_2}. We used formula 6.8(26) in
\cite{Bateman_tab_vol_1} ($\mathrm{Re}\,c > |\mathrm{Re}\,\nu|$)
\begin{equation}\label{contour_int_A3}
\frac{1}{2\pi i}\int\limits_{c - i\infty}^{c + i\infty} ds \,
2^{\,s-2}\,\Gamma \!\left( \frac{s+\nu}{2}\right)\Gamma \!\left(
\frac{s-\nu}{2}\right) \,(ax)^{-s} = K_\nu(ax) \;.
\end{equation}
\textbf{A.4.} $f(x) = 1/(x^2+a^2)^{n+3/2}$, $a>0$, integer
$n\geqslant0$. Then we get the integral
\begin{equation}\label{example_int_A4}
A_4(q) = \int\limits_0^{\infty} \!\!\frac{x}{(x^2 + a^2)^{n+3/2}}
\,J_0(q x) \,dx \;.
\end{equation}
The derivatives of $g(z)$ are: $\bar{g}^{(m)}(0) = a^{-2m-2n-3}
\Gamma(m+n+3/2)/\Gamma(n+3/2)$. We obtain from \eqref{A_cont_final}
($-1 < \alpha <0$)
\begin{equation}\label{A4_1}
A_4(q) = \frac{1}{\Gamma(n+3/2)} \,\frac{1}{a^{2n+3} q^2}
\frac{1}{\pi i} \int\limits_{\alpha - i\infty}^{\alpha + i\infty}
\!\!ds  \,\Gamma(s+1)\Gamma(s+n+3/2) \left(\frac{q a}{2}
\right)^{\!\!-2s} \;.
\end{equation}
By using formula \eqref{contour_int_A3}, after change of variable
$s+1= (s'/2-n-1/2)/2$, we find that
\begin{equation}\label{A4_2}
A_4(q) = \frac{1}{\Gamma(n+3/2)} \left(
\frac{q}{2a}\right)^{\!n+1/2} \!\!K_{n+1/2}(q a) \;,
\end{equation}
see eq.~2.12.4.28. in \cite{Prudnikov_vol_2}.\\
\textbf{A.5.} $f(x) = J_0(a x)/(x^2+c^2)$, where $q>a>0$, $c>0$.
Then we come to the integral
\begin{equation}\label{example_int_A5}
A_5(q) = \int\limits_0^{\infty} \!\!\frac{x}{x^2 + c^2} \,J_0(a x)
\,J_0(q x) \,dx  \;.
\end{equation}
The derivatives of $g(z)$ are the following:
\begin{align}\label{der_example_int_A5}
\bar{g}^{(m)}(0) &=c^{-2m-2} \,\Gamma(m+1) \sum_{k=0}^m
\frac{1}{(k!)^2} \left( \frac{ac}{2} \right)^{\!2k} \nonumber \\
&= c^{-2m-2} \,\Gamma(m+1) I_0(ac)
\nonumber \\
&- \frac{1}{(m+1) \Gamma(m+2)} \,{}_1F_2 \!\left(1; m+2,m+2;
\frac{a^2c^2}{4} \right) \left( \frac{a}{2} \right)^{\!2m+2}.
\end{align}
Here and in what follows ${}_pF_q(a_1, \ldots, a_p; b_1, \ldots,
b_q; z)$ is the generalized hypergeometric series
\cite{Bateman_vol_1}. Then we obtain from \eqref{A_cont_final},
\eqref{der_example_int_A5}
\begin{equation}\label{example_2.5_1}
A_5(q) = A_5^{(1)}(q) + A_5^{(2)}(q) \;,
\end{equation}
where (see eq.~7.3(17) in \cite{Bateman_tab_vol_1})
\begin{equation}\label{example_A5_A5_1}
A_5^{(1)}(q) = \frac{I_0(ac)}{q^2 c^2} \frac{1}{\pi i}
\int\limits_{\alpha - i\infty}^{\alpha + i\infty} \!\!ds \,
\Gamma^2(s+1) \left( \frac{q^2 c^2}{4} \right)^{\!\!-s} \! = I_0(ac)
K_0(q c) \;,
\end{equation}
and
\begin{equation}\label{example_A5_A5_2}
A_5^{(2)}(q) =  -\frac{1}{\pi i q^2} \int\limits_{\alpha -
i\infty}^{\alpha + i\infty} \!\!ds \, \frac{1}{(s+1)^2} \,{}_1F_2
\!\left(1; s+2,s+2; \frac{a^2c^2}{4} \right) \!\left(
\frac{q^2}{a^2} \right)^{\!\!-s} \!,
\end{equation}
with $-1 < \alpha <0$. The integrand in \eqref{example_A5_A5_2} is a
regular function in the half-plane $\mathrm{Re}s>-1$. Moreover,
since $q>a$, it decreases very rapidly as $|s|\rightarrow \infty$ in
the half-plane $\mathrm{Re}s>-1$. According to the Cauchy integral
theorem, $A_5^{(2)}(q) = 0$.

As a result,
\begin{equation}\label{A5}
A_5(q) = I_0(ac) K_0(q c) \;,
\end{equation}
see also eq.~2.12.32.11. in \cite{Prudnikov_vol_2}.

Our formula \eqref{A_cont_final} enables one to calculate new
improper parametric integrals of Bessel functions which are not yet
presented in the most complete tables of integrals
\cite{Prudnikov_vol_2}, \cite{Gradshteyn}. \\
\textbf{A.6.} As an example, let us consider the following integral
($a, c>0$)
\begin{equation}\label{example_int_A6}
A_6(q) = \int\limits_0^{\infty} \!\!\frac{x}{x^2 + c^2} \, e^{-a^2
x^2} \!J_0(q x) \,dx \;,
\end{equation}
i.e. zero-order Hankel transform with the function $f(x) = e^{-a^2
x^2}/(x^2 + c^2)$. We find that
\begin{equation}\label{der_example_int_A6}
\bar{g}^{(m)}(0) = c^{-2m-2} \,\Gamma(m+1) \sum_{k=0}^m \frac{(a
c)^{2k}}{k\,!} \;.
\end{equation}
Using relation
\begin{equation}\label{sum_trans}
\sum_{k=0}^m \frac{(a c)^{2k}}{k\,!} = e^{(a c)^2} - (ac)^{2(m+1)}
\sum_{p=0}^\infty \frac{(a c)^{2p}}{\Gamma(p+m+2)} \;,
\end{equation}
we obtain
\begin{equation}\label{I6_two_terms}
A_6(q) = A_6^{(1)}(q) + A_6^{(2)}(q) \;,
\end{equation}
where (see eq.~\eqref{example_A5_A5_1})
\begin{equation}\label{A6_1}
A_6^{(1)}(q) = \frac{e^{(a c)^2}}{\pi i(q c)^2} \int\limits_{\alpha
- i\infty}^{\alpha + i\infty} \!\!ds \, \Gamma^2(s+1) \left(
\frac{q^2 c^2}{4} \right)^{\!\!-s} \! = e^{(a c)^2} K_0(q c) \;,
\end{equation}
and
\begin{equation}\label{A6_2_1}
A_6^{(2)}(q) = -\frac{a^2}{\pi iq^2} \sum_{p=0}^\infty \,(a c)^{2p}
\int\limits_{\alpha - i\infty}^{\alpha + i\infty} \!\!ds \,
\frac{\Gamma^2(s+1)}{\Gamma(s + p +2)} \left( \frac{q^2}{4a^2}
\right)^{\!\!-s} ,
\end{equation}
with $-1 < \alpha <0$. According to eqs.~7.3(43) from
\cite{Bateman_tab_vol_1} and 5.6(6) from \cite{Bateman_vol_1}, we
have
\begin{equation}\label{A6_2_2}
\frac{1}{2\pi i}\int\limits_{c - i\infty}^{c + i\infty} \!\!ds \,
\frac{\Gamma^2(s+1)}{\Gamma(s + p +2)} \left( \frac{q^2}{4a^2}
\right)^{\!\!-s} \! = \frac{q}{2a} \,e^{-q^2/8 a^2} W_{-p - 1/2, \,
0} \!\left( \frac{q^2 }{4a^2} \right) ,
\end{equation}
where $W_{\chi, \mu}(x)$ is the Whittaker function of the second
kind \cite{Bateman_vol_1}. We can apply the relation 6.9(5) from
\cite{Bateman_vol_1}:
\begin{equation}\label{Wittaker_vs_hypergeom_1}
W_{\chi, \mu}(x) = e^{-x/2} x ^\chi \,{}_2F_0 \!\left( \frac{1}{2}
-\chi + \mu, \frac{1}{2} -\chi - \mu; -\frac{1}{x} \right) .
\end{equation}
After that we find the analytic expression for
\eqref{example_int_A6}
\begin{align}\label{A6_2_3}
A_6(q) &= e^{(a c)^2} K_0(q c) - \frac{2 a^2}{q^2} \,e^{-q^2/4a^2}
\nonumber \\
&\times \sum_{p=0}^\infty \left( \frac{2 a^2c}{q} \right)^{\!\!2p}
\!{}_2F_0 \!\left(p+1, p+1; -\frac{4a^2}{q^2} \right) .
\end{align}

For $a=0$ it coincides with the known integral (see eq.~2.12.4.23 in
\cite{Prudnikov_vol_2})
\begin{equation}\label{example_2.6_a_zero}
A_6(q)\Big|_{a = 0} = \int\limits_0^{\infty} \!\!\frac{x}{x^2 + c^2}
\, J_0(q x) \,dx = K_0(q c) \;.
\end{equation}
In the limit $c \rightarrow 0$ the leading term of
\eqref{example_int_A6},
\begin{equation}\label{A6_small_c}
A_6(q)\Big|_{c \rightarrow 0} \simeq \int\limits_0^{\infty}
\!\!\frac{y}{y^2 + 1} \, J_0(q c y) \,dy  = K_0(q c) \;,
\end{equation}
is also reproduced by eq.~\eqref{A6_2_3}.

Let us demonstrate that our formula \eqref{A6_2_3} gives a correct
result for $q=0$. To do this, we use the following equation
\cite{Bateman_vol_1}
\begin{equation}\label{Wittaker_vs_hypergeom_2}
W_{\chi, \mu}(x) = e^{-x/2} x ^{\mu+1/2} \,\Psi \!\left( \frac{1}{2}
-\chi + \mu , 2\mu + 1; x \right) ,
\end{equation}
where $\Psi(a,c,x)$ is the confluent hypergeometric series of the
second kind \cite{Bateman_vol_1}. Then we get another representation
for \eqref{example_int_A6}
\begin{equation}\label{A6_2_4}
A_6(q) = e^{(a c)^2} K_0(q c) - \frac{1}{2} \,e^{-q^2/4a^2}
\sum_{p=0}^\infty (a c)^{2p} \Psi\!\left(p+1, 1; \frac{q^2}{4a^2}
\right) .
\end{equation}
The function $\Psi(a,1,x)$ has the following asymptotic
behavior as $x\rightarrow 0$ \cite{Abramowitz}:%
\footnote{In the similar formula 6.8(5) in \cite{Bateman_vol_1} the
term $2 \gamma_E$ has the wrong sign.}
\begin{equation}\label{confl_small_x_2.7}
\Psi(a,1,x) = - \frac{1}{\Gamma(a)} \,[\,\ln x + \psi(1) + 2
\gamma_E] + \mathrm{O}(x \ln x) \;,
\end{equation}
where $\psi(x)$ is the psi function, and $\gamma_E$ is the
Euler--Mascheroni constant. As a result, we get in the limit $q
\rightarrow 0$
\begin{align}\label{series_small_q_A6}
&-\frac{1}{2} \sum_{p=0}^\infty (a c)^{2p} \Psi\!\left(p+1, 1;
\frac{q^2}{4a^2} \right)
\nonumber\\
&= e^{(a c)^2} [\, \ln (q/2a) + \gamma_E ] + \frac{1}{2}
\sum_{p=0}^\infty (a c)^{2p} \frac{\psi(p+1)}{p\,!} \;.
\end{align}
It is known that
\begin{equation}\label{sum_psi}
\sum_{p=0}^\infty \frac{\psi(p+1)}{\Gamma(p+1)} \, x^{p} = e^{x} \,
[\Gamma (0, x) + \ln x] \;,
\end{equation}
where $\Gamma (a, x)$ is the incomplete gamma function
\cite{Bateman_vol_2}. Since
\begin{equation}\label{K0_small_q}
K_0(x) = - \gamma_E - \ln(x/2) + \mathrm{O}(x^2) \;,
\end{equation}
we obtain from eqs.~\eqref{A6_2_4}, \eqref{series_small_q_A6},
\eqref{sum_psi}:
\begin{equation}\label{A6_q_zero}
A_6(0) = \frac{1}{2}  \, e^{(a c)^2} \, \Gamma \left( 0, (ac)^2
\right) \;,
\end{equation}
in full accordance with eq.~2.3.4.3 from \cite{Prudnikov_vol_1}
\begin{equation}\label{example_int_A6_q_zero}
A_6(q)\Big|_{q = 0} = \int\limits_0^{\infty} \!\!\frac{x}{x^2 + c^2}
\, e^{-a^2 x^2} dx = \frac{1}{2} \, e^{(a c)^2} \, \Gamma \!\left(
0, (ac)^2 \right) \;.
\end{equation}

For large $a$ or $c$ we have to take $\bar{g}^{(m)}(0) = c^{-2}
a^{2m}$, then
\begin{equation}\label{A6_large_a}
A_6(q)\big|_{a \gg 1} \simeq \frac{1}{\pi i(q c)^2}
\int\limits_{\alpha - i\infty}^{\alpha + i\infty} \!\!ds \,
\Gamma(s+1) \left( \frac{q^2}{4a^2} \right)^{\!\!-s} \! =
\frac{1}{2(a c)^2} \,e^{-q^2/4 a^2} \;,
\end{equation}
where $-1 < \alpha <0$. On the other hand,
\begin{equation}\label{example_int_2.6_q_a_large}
A_6(q) \Big|_{ac \gg 1} \simeq \frac{1}{(a c)^2}
\int\limits_0^{\infty} \!\!z \, e^{- x^2} \!J_0(qz/a) \,dz =
\frac{1}{2(a c)^2} \,e^{-q^2/4 a^2} \;.
\end{equation}

Finally, we predict the following asymptotic behavior
\begin{equation}\label{A6_large_q}
A_6(q)\big|_{q \gg 1} \stackrel{\mathrm{as}}{=} \, e^{(a c)^2} K_0(q
c) \;.
\end{equation}

Let us underline that integrals $A_i(q)$ ($i=1,2, \ldots 6$) are
exponentially decreasing functions at large $q$. It is in agreement
with the suggestion made in \cite{Frenzen} about the integral
\begin{equation}\label{Frenzen_int}
I_f(q) = \int\limits_0^\infty \!\!x \,e^{-x^2} \!\varphi(x^2) J_0(q
x) \,dx \;.
\end{equation}
It says that if functions $\varphi(z)$ are entire or meromorphic,
$I_f(q) = \mathrm{O}(e^{-\gamma q^2})$ and $I_f(q) =
\mathrm{O}(e^{-\delta q})$, respectively, as $q\rightarrow \infty$.
See asymptotics of our integrals in the examples A.1-A.2 and
A.3-A.6, correspondingly.

\textsc{Theorem 2.} The Hankel transform \eqref{A} with the function
$f(x) = h(x^4)$ and $q>0$ may be expressed in the form
\begin{equation}\label{A_cont_2_final}
A(q) = \frac{1}{q^2 \sqrt{\pi} i} \!\int\limits_{\alpha -
i\infty}^{\alpha + i\infty} \!\! ds \, \bar{h}(s) \,
\frac{\Gamma(2s+1)}{\Gamma(1/2-s)} \,2^{\,6s+1} \,q^{-4s} \;,
\end{equation}
for $-1/2 < \alpha < 0$, provided that

\begin{description}
  \item[(1)]
$h(z)$ is a regular function and its Taylor series at $z=0$  has the
form
\begin{equation}\label{h_Taylor_series}
h(z) = \sum_{m=0}^\infty \frac{\bar{h}^{(m)}(0)}{m!} \,(-z)^m \;;
\end{equation}
  \item[(2)]
$h(z)= \mathrm{O}(z^{-d})$, as $z\rightarrow \infty$, for $d>1/8$;
  \item[(3)]
$\bar{h}^{(s)}(0)$ is a regular (single-valued) function defined on
a half-plane
\begin{equation}\label{h_half_plane}
H(\delta) = \{s\in \mathrm{\mathbb{C}}\!: \mathrm{Re} \,s \geqslant
-\delta \}
\end{equation}
for some $1/8 < \delta < 1/2$ and satisfies the growth condition
\begin{equation}\label{h_growth}
|\bar{g}^{(s)}(0)| < C e^{P v + A|w|}
\end{equation}
for some $A < \pi/2$ and all $s=v+iw \in H(\delta)$.
\end{description}

\textsc{Proof}. The proof proceeds as the proof of the Theorem~1.

Let us consider one example.\\
\textbf{A.7.} $f(x) = 1/\sqrt{x^4+a^4}$, $a>0$. Then we have the
following integral
\begin{equation}\label{example_int_A7}
A_7(q) = \int\limits_0^{\infty} \!\!\frac{x \,J_0(q x)}{\sqrt{x^4 +
a^4}} \,dx \;.
\end{equation}
Since $\bar{h}^{(m)}(0) = (\pi)^{-1/2} a^{-4m-2} \,\Gamma(m + 1/2)$,
we find from \eqref{A_cont_2_final} that
\begin{equation}\label{A7_1}
A_7(q) = \frac{1}{(q a)^2} \frac{1}{\pi i}\int\limits_{\alpha -
i\infty}^{\alpha + i\infty} \!ds \, 2^{\,6s+1} \,\frac{\Gamma(2s+1)
\Gamma (s + 1/2)}{\Gamma(1/2-s)} \, (q a)^{-4s} \;.
\end{equation}
Let us put $s = s'/4-1/2$, then ($0 < \alpha' < 2$)
\begin{equation}\label{A7_2}
A_7(q) = \frac{1}{2\pi i}\int\limits_{\alpha' - i\infty}^{\alpha' +
i\infty} \!ds' \, 2^{\,s'-3} \,\frac{\Gamma(s'/2) \Gamma
(s'/4)}{\Gamma(1-s'/4)} \, \left( \frac{q a}{\sqrt{2}}
\right)^{\!\!-s'}.
\end{equation}
Now we can apply formula 6.8(41) from ref.~\cite{Bateman_tab_vol_1}
to find that
\begin{equation}\label{A7_3}
A_7(q)  = J_0\!\left( \frac{q a}{\sqrt{2}} \right) \!K_0\!\left(
\frac{q a}{\sqrt{2}} \right) ,
\end{equation}
in accordance with eq.~2.12.5.2. in \cite{Prudnikov_vol_2}. Again,
$A_7(q)$ is an exponentially decreasing function at large $q$, up to
damped oscillations of $J_0(q a/\!\sqrt{2})$.



\section{Conclusions} %

With the use of Ramanujan's Master Theorem, we have derived two
formulas which define the Hankel transforms of order zero
$\mathcal{H}_0(q)$  \eqref{A} with the even functions $g(x^2)$,
$h(x^4)$ in terms of the inverse Mellin transforms
\eqref{A_cont_final}, \eqref{A_cont_2_final}, provided that
functions $g(z)$, $h(z)$ and their derivatives satisfy special
conditions. The representation obtained is valid not only for $q
\rightarrow \infty$ but for all $q>q_0>0$, if these conditions are
satisfied.%
\footnote{A particular value of $q_0$ depends on the form of $g(z)$
or $h(z)$, see examples A.1-A.7.}
The validity of our formulas is illustrated by a number of examples.
It is shown that formula \eqref{A_cont_final} can be useful in
calculating as yet unknown parametric improper integrals of the
Bessel function $J_0(x)$ (see integral \eqref{example_int_A6} and
eqs.~\eqref{A6_2_3}, \eqref{A6_2_4}).



\section*{Acknowledgements}

The author is indebted to A.P.~Samokhin for fruitful discussions.



\setcounter{equation}{0}
\renewcommand{\theequation}{A.\arabic{equation}}

\section*{Appendix A}
\label{app:A}

Let us calculate the derivatives of the function ($n = 1, 2,
\ldots$),
\begin{equation}\label{g}
g(z) = e^{-a^2 z} J_{2n}(c \sqrt{z}) \;,
\end{equation}
taken at $z=0$. We start from the expression
\begin{equation}\label{der_g2}
g^{(m)}(z) = \sum_{p=0}^m \binom{m}{p} \!\left[ J_{2n}(c \sqrt{z})
\right]^{(p)} \left[ e^{-a^2 z} \right]^{(m-p)} ,
\end{equation}
where
\begin{equation}\label{der_Bessel2n_zero}
\left[ J_{2n}(c \sqrt{z}) \right]^{(p)} \Big|_{z=0} = \frac{(-1)^{p}
\,p\,!}{(p-n)!\,(p+n)! } \left( \frac{c^2}{4} \right)^{\!\!p} .
\end{equation}
As a result, we obtain that $g^{(m)}(0)=0$ for $0 \leqslant m
\leqslant n-1$, while for $m \geqslant n$
\begin{align}\label{der_g2_final}
g^{(m)}(0) &= (-1)^{m} a^{2m} \sum_{p=n}^m \binom{m}{p}
\frac{p\,!}{(p-n)!\,(p+n)!} \left( \frac{c^2}{4a^2} \right)^{\!\!p}
\nonumber \\
&= (-1)^m \frac{m!}{(m+n)!} \,a^{2m} \!\left( -\frac{c^2}{4 a^2}
\right)^{\!\!n} \!L_{m-n}^{2n} \!\left( - \frac{c^2}{4 a^2} \right)
,
\end{align}
where
\begin{equation}\label{laguerre_gen}
L_n^\alpha (z) = \sum_{k=0}^m \binom{n+\alpha}{n-k}
\frac{(-z)^k}{k\,!}
\end{equation}
is the generalized Laguerre polynomial, $\alpha > -1$
\cite{Bateman_vol_2}. In particular, we find for $n=0$
\begin{equation}\label{der_g3_final}
g^{(m)}(0) = (-1)^m \,a^{2m} L_{m} \!\left( - \frac{c^2}{4 a^2}
\right),
\end{equation}
where $L_n(x) = L_n^0(x)$.





\begin{thebibliography}{99}
%
\bibitem{Tranter}
C.J.~Tranter, \emph{Integral Transforms in Mathematical Physics},
John Wiley \& Sons, New York, 1951.
\bibitem{Sneddon}
Ian N.~Sneddon, \emph{Fourier Transforms}, Dover Publications, New
York, 1995 (first published 1951).
\bibitem{Debnath}
L.~Debnath and D.~Bhatta, \emph{Integral Transforms and Their
Applications}, Third Edition, CRC Press, New York, 2015.
%
\bibitem{Moliere}
G.~Moli\`{e}re, Z. Naturforsch, \textbf{2a}, (1947) 133; Z.
Naturforsch, \textbf{3a} (1948) 78.
\bibitem{Glauber}
R. J. Glauber, in \emph{Lectures in Theoretical Physics}, ed.
W.E.~Brittin and L.G.~Dunham, New York, Vol.~1, 1959, p.~315.
\bibitem{Collins}
P.D.B.~Collins, \emph{An introduction to Regge theory} \& \emph{high
energy physics}, Cambridge University Press, 1977.
%
\bibitem{Willis}
H.F.~Willis, Phil. Mag. \textbf{39} (1948) 455.
%
\bibitem{Hsu}
L.C.~Hsu, Ann. Pol. Math. \textbf{11} (1961) 7.
\bibitem{Slonovskii}
N.V.~Slonovskii, Izv. Vyssh. Uchebn. Zaved. Mat. \textbf{5} (1968)
86.
\bibitem{Handelsman}
R.A.~Handelsman and J.S.~Lew, J. Math. Anal. Applic. \textbf{35}
(1971) 405.
\bibitem{MacKinnon}
R.F.~Mackinnon, Math. Comp. \textbf{26} (1972) 515.
\bibitem{Wong}
R.~Wong, SIAM J. Math. Anal. \textbf{7} (1976) 799.
\bibitem{Soni}
K.~Soni, \emph{Asymptotic Expansion of the Hankel Transform with
Explicit Remainder Terms}, Quart. Appl. Math., \textbf{40} (1982) 1.
%
\bibitem{Frenzen}
C.L.~Frenzen and R.~Wong, Math. Comp. \textbf{45} (1985) 537.
%
\bibitem{Hardy}
G.H.~Hardy, \emph{Ramanujan. Twelve Lectures on Subjects Suggested
by His Life and Work}, 3nd edn., Celsea, New York, 1978.
\bibitem{Amdeberham}
T.~Amdeberham, J.~Espinosa, I.~Gonzalez \emph{et al.}, The Ramanujan
J. \textbf{29} (2012) 103.
%
\bibitem{Bateman_vol_2}
\emph{Higher Transcendental Functions}. Vol.~2. By the staff of the
Bateman manuscript project (A.~Erd\'{e}lyi,\emph{ Editor};
W.~Magnus, F.~Oberhettinger, F.G.~Tricomi, \emph{Associates}),
McGraw-Hill Book Company, New York, 1953.
%
\bibitem{Bateman_vol_1}
\emph{Higher Transcendental Functions}. Vol.~1. By the staff of the
Bateman manuscript project (A.~Erd\'{e}lyi,\emph{ Editor};
W.~Magnus, F.~Oberhettinger, F.G.~Tricomi, \emph{Associates}),
McGraw-Hill Book Company, New York, 1953.
%
\bibitem{Prudnikov_vol_2}
A.P. Prudnikov, Yu.A. Brychkov and O.I. Marichev, \emph{Integrals
and Series}, Vol. 2: \emph{Special Functions}, Gordon \& Breach Sci.
Publ., New York, 1986.
%
\bibitem{Bateman_tab_vol_1}
\emph{Tables of Integral Transforms}. Vol.~1. By the staff of the
Bateman manuscript project (A.~Erd\'{e}lyi,\emph{ Editor};
W.~Magnus, F.~Oberhettinger, F.G.~Tricomi, \emph{Associates}),
McGraw-Hill Book Company,  New York, 1954.
%
\bibitem{Gradshteyn}
I.S.~Gradshteyn and I.M.~Ryzhik, \emph{Table of Integrals, Series,
and Products}, 7th edition, Academic Press, New York, 2007.
%
\bibitem{Abramowitz}
M.~Abramowitz and I.A.~Stegun (Eds.), \emph{Handbook of Mathematical
Functions with Formulas, Graphs, and Mathematical Tables}, Tenth
printing, Dover, New York, 1972.
%
\bibitem{Prudnikov_vol_1}
A.P. Prudnikov, Yu.A. Brychkov and O.I. Marichev, \emph{Integrals
and Series}, Vol. 1: \emph{Elementary Functions}, Gordon \& Breach
Sci. Publ., New York, 1986.
%
\end{thebibliography}
\end{document}